# A calculus for causal inference with instrumental variables

Wing Hung Wong

Department of Statistics and Department of Biomedical Data Science

Stanford University, Stanford, CA 94305

**Abstract**:

Under a general structural equation framework for causal inference, we provide a definition of the causal effect of a variable X on another variable Y, and propose an approach to estimate this causal effect via the use of instrumental variables.

## 1) Introduction

The purpose of this paper is to discuss how to infer the causal effect of an explanatory variable $X$ on a response variable $Y$. It is assumed that $Y$ is related to $X$ by a structural equation $Y = f(X, U)$ where, $Y, X, U$ are regarded as random variables in some probability space $(\Omega, A, P)$. Here $U$ represents all other variables besides $X$ that may affect $Y$. Our objective to define the causal effect of $X$ on $Y$ in this context, and to develop a method to estimate this effect with the help of an instrumental variable $Z$. We begin by illustrating the issues in some examples.

Example 1: Suppose $\tilde{Y} = \beta\tilde{X}$, $X = \tilde{X} + U_1$, $Y = \tilde{Y} + U_2$, where $U_1, U_2$ are zero-mean errors uncorrelated with $\tilde{X}$ and with each other. Given observations on $X$ and $Y$, how to estimate the parameter $\beta$ which is regarded as the causal effect of $\tilde{X}$ on $\tilde{Y}$? This is the classic problem of errors in variables (reviewed in Durbin 1954). We can write

$$Y = f(X, U) = \beta X + (U_2 - \beta U_1)$$

Since the error term $U_2 - \beta U_1$ affects the distributions of both $X$ and $Y$, the regression $Y$ on $X$ does not provide a consistent estimate for $\beta$. However, suppose we can find a variable $Z$ such that

$$X = \alpha Z + U_3$$

where Z is uncorrelated with $U_1, U_2, U_3$ but $U_3$ may be correlated with $U_1, U_2$. Then, Z can be used as an "instrumental variable" to construct the estimate $\hat{\beta} = b/a$, where $b$ is obtained by the linear regression of $Y$ on $Z$ based on observed data on $(Y, Z)$, and $a$ is similarly obtained from the regression of $X$ on $Z$. It is easy to show that $\hat{\beta}$ is a consistent estimate of $\beta$.

Example 2: Consider the following model for econometric program evaluation (Heckman and Robb 1985)

$$Y = \beta_0 + \beta_1 X + U_1, \quad X = I(X^* > 0), \quad X^* = \alpha_0 + \alpha_1 Z + U_2$$

where $X$ is the indicator of whether a subject participated in the program under study and $Y$ represents outcome for the subject. The subject's decision on whether to participate depends on whether the latent variable $X^*$ exceeds a threshold, and $X^*$ itself depends on the variable $Z$, which is assumed to be independent of the errors $U_1$ and $U_2$. The errors are allowed to be dependent, so $X$ is not uncorrelated with $U_1$. In this model, we are interested in $\beta_1$, which represents the change in outcome due participation to the program. With correlated $X$ and $U_1$ this causal effect cannot be revealed directly by fitting $E(Y|X)$. However, the authors pointed out that a consistent estimate of $\beta_1$ can be derived by treating $Z$ as an instrumental variable. We will return to this example in section 3.

Example 3: In the Neyman-Rubin potential outcome framework (Rubin 1974), a randomized experiment to compare the effects of two treatments can be analyzed with the help two variables $Y_i(x), x \in \{0,1\}$ that represent the potential outcomes of subject $i$ under treatment 0 and treatment 1 respectively. Then, the population-level causal effect $\theta$ is defined as the population average of the subject-level treatment effect $d_i = Y_i(1) - Y_i(0)$. Note that $d_i$ is a counterfactual variable that is not observable for any subject but is needed to define the populational-level causal effect. For simplicity, hereafter the term causal effect will refer to population-level causal effect, while subject-level causal effect will be written explicitly as such. Under the assumption that the treatment $X_i$ for subject $i$ is chosen randomly from the two available treatments, one can estimate $\theta$ by the difference of the two treatment group means, i.e. $\hat{\theta} = \bar{Y}_1 - \bar{Y}_2$. In practice, however, the actual treatment $X_i$ received by subject $i$ may differ from the (randomly) assigned treatment $Z_i$. In this case, $\hat{\theta}$ is not a suitable estimate of the treatment effect, but one can try to use $Z$ as an instrumental variable to construct a better estimate. It turns out that $\theta$ is not identifiable from observations on $(Y, X, Z)$ without additional conditions. Previous authors had formulated conditions under which one can use $Z$ as an instrumental variable to obtain estimates of various "treatment effects" that represent averages of $d_i$ over suitably defined subpopulations (Heckman 1990, Imbens and Angrist 1994, Angrist, Imbens and Rubin 1996). We can analyze this example under the structural equation framework by thinking of $\Omega$ as a large population of subjects, and replace the subject index $i$ by a sample point $\omega \in \Omega$. Let $Z(\omega)$ and $X(\omega)$ be respectively the assigned treatment and the actual treatment for subject $\omega$, and consider the structural equation $Y = f(X, U) = XU_1 + (1-X)U_0$. Then, the potential outcome for subject $\omega$ under treatment 1 is given by $f(1, U(\omega)) = U_1(\omega)$. Similarly, the potential outcome for subject $\omega$ under treatment 0 is given by $U_0(\omega)$. Therefore $\theta = E(U_1 - U_0)$ is the causal effect of interest. We will return to this example in section 2.

Example 4: In this example, $Y = s(X) + U_1$, $X = t(Z) + U_2$, where $X, Z, U_1, U_2$ are real-valued random variables. We assume that $s$ and $t$ are smooth functions and $Z$ is independent of $U_1, U_2$, but the distribution of $(U_1, U_2)$ is allowed to be arbitrary. This is an extension of models with linear equations such as example 1. In this situation, what is the causal effect of $X$ on $Y$, and how to estimate it using instrumental variables? We note that in modern applications it is often true that the sample size is large enough to infer functions nonparametrically. Thus a causal inference methodology to accommodate the use of non-linear, nonparametric models will be valuable. We will return to this example in Section 4.

In this paper, we will extend the method of instrumental variable to cover a larger class of models than have been treated previously. Both $X$ and $Z$ are allowed to be multidimensional and the equations $Y = f(X, U)$ and $X = g(Z, U)$ can be nonlinear. Throughout, we assume that $Y$ is either real-valued or binary. Extension to discrete $Y$ is possible but does not seem to provide additional insight. Although the potential outcome framework has also been widely used for causal inference, we use the structural equation framework to develop our theory as it is easier to work with continuous variables and nonlinear models under this framework. As seen in example 3 above and further discussed in section 2, the two frameworks are largely equivalent conceptually. Based on the structural equation $f( )$, we define the causal effect of $X$ on $Y$ as the expected change of $Y$ induced by changing $X$ when its current value $x$, while keeping $U$ fixed. Denoting this causal effect by $\theta(x)$, we develop conditions under which $Z$ can be used as an instrumental variable to construct consistent estimates of the function $\theta$. One of the key condition can be stated as "the change of $Y$ caused by varying $X$ should be uncorrelated with the change of $X$ caused by varying $Z$" (or conditionally uncorrelated given Z, in the case of continuous $X$). Under these conditions, we explain how to obtain $\theta$ in terms of $E(Y|Z=z)$ and $E(X|Z=z)$. In this way, the method of instrumental variable is extended to cover a large class of nonlinear, nonparametric causal models.

To our knowledge these results are new. However there are a number of prior works on nonlinear, nonparametric models that are highly relevant. Newey and Powell (2003) studied the model $Y = f(X, Z_1) + U$, where $Z_1$ is scalar, $E(U|Z) = 0$, $Z = (Z_1, Z_2)$. They showed that $f(\ )$ is identifiable from the distribution of $(Y, X, Z)$ if the family of conditional distributions $X$ given $Z$ is "complete" in the sense that $E(r(X, z_1)|Z = (z_1, z_2)) = 0$ for all $z_2$ implies that $r(X, z_1) = 0$ in distribution. Chernozhukov, Imbens and Newey (2007) considered the model $Y = f(W, U)$, $W = (X, Z_1), Z = (Z_1, Z_2)$ where $U$ is a scalar latent variable, $f$ is required to be monotone in $U$, and $Z$ is independent of $U$. Imbens and Newey (2009) studied the model $Y = f(X, \varepsilon)$, $X = (X_1, Z_1), X_1 = h(Z_1, Z_2, U)$, where $U$ is a scalar disturbance and $h$ is required to be monotone in $U$. Their interest was to construct a "control variable" $V$ based on the conditional distribution $X_1|Z$. Once $V$ is available, certain structural features (different from $\theta$) can be identified in terms of the conditional distribution $Y|(X, V)$. Torgovitsky (2019) studied the model $Y = f(X, U_0), X = (X_1, .., X_p), X_k = g_k(Z, U_k), k = 1, ... p$; where $U_0, U_1, ... U_p$ are scalar latent variables, $f$ is required to be monotone in $U_0$, and $g_k$ is required to be monotone in $U_k$. Finally, in the case of binary $X$, Kennedy, Lorch and Small (2018) studied a nonlinear and nonparametric causal model, and developed methods to estimate the local instrumental variable curve. Although these papers contained identifiability results, the features they identified are different from the causal effect $\theta(x)$ of interest here. Also, these papers made use of various monotonicity assumptions on the dependency on the $U_k's$, which may limit the applicability of the results.

The plan for the reminder of the paper is as follows. In section 2 we discuss the case when both $X$ and $Z$ are discrete. In section 3 we discuss the case when $X$ is discrete and $Z$ is continuous. In section 4 we discuss the case when both $X$ and $Z$ are continuous. In section 5, we provide the proofs of the theorems.

### 2) The case of discrete $X$ and $Z$

Suppose Z takes value in $\{z_1, ..z_m\}$ and X takes value in $\{x_0, x_1, ..x_n\}$, and Y, X, Z are related by the following structural equations

$$X = \sum_1^m I(Z = z_i)V_i \quad , Y = \sum_0^n I(X = x_j)U_j \qquad (1)$$

Similar to example 3, we can interpret this model in the usual potential outcome framework of causal inference by regarding $U_j(\omega)$ as the value of $Y$ for subject $\omega$ when $X(\omega)$ takes the value $x_j$, Then $U_j(\omega) - U_k(\omega)$ is the change in $Y$ for subject $\omega$ when $X$ changes from $x_k$ to $x_j$, and we can regard the population average of these changes, $E(U_j - U_k)$, as the causal effect on $Y$ of when $X$ changes from $x_k$ to $x_j$. Suppose $x_0$ is regarded as the baseline state for $X$, then the causal effects are determined by the vector $\theta = (\theta_1, ... \theta_n)$ where $\theta_j = E(U_j - U_0)$. At this stage we do not impose any conditions on the joint distributions of the $U's$ and $V's$. Thus (1) provides a general model for the dependency of $Y$ on $X$ as well as the dependency of $X$ on $Z$.

Because $U_j$ and $X$ are not independent variables, $E(Y|X = x_j) = E(U_j|X = x_j) \neq E(U_j)$, and it is not possible to estimate $\theta$ based on observation only on $X$ and $Y$. However, under some conditions, $\theta$ is identifiable by $(X, Y, Z)$, i.e. they are determined by the joint distribution of $(X, Y, Z)$ and can be estimated consistently from observations on these variable. To formulate these conditions, let

$\Delta_{j0}Y = U_j - U_0 = $ change in $Y$ when $X$ changes from $x_0$ to $x_j$

$I_j(X) = I(X = x_j)$

$\Delta_{ik}I_j(X) = I(V_i = x_j) - I(V_k = x_j) = $ change in $I(X = x_j)$ when Z changes from $z_k$ to $z_i$.

Let $S$ be the set of unordered pairs of indexes from $\{1, ..., m\}$. For any $s = \{i, k\} \in S$, define

$b_s = E(Y|Z = z_i) - E(Y|Z = z_k)$

$a_{sj} = P(X = x_j|Z = z_i) - P(X = x_j|Z = z_k)$.

<u>Theorem 1</u>: Suppose Z is independent of $\{V_1, ... V_m; U_0, U_1, ... U_n\}$, and there exists $S_0 \subset S$ such

i)      $\Delta_{j0}Y \perp \Delta_{ik}I_j(X)$ for all $\{i, k\} \in S_0$ and $j \in \{1, ..n\}$, where $\perp$ means zero correlation
ii)      $Rank(A) \geq n$, where $A$ is the matrix $\{a_{sj} : s \in S_0, j = 1,..n\}$

then $\theta$ is identifiable via the equation $A\theta = b$ where $b$ is the vector $\{b_s : s \in S_0\}$

Since $A$ and $b$ can be obtained from the conditional distribution of $X|Z$ and $Y|Z$ respectively, $Z$ can thus be used as an instrumental variable to estimate the causal effect of $X$ on $Y$. Note that in addition to being independent of $\{V_1, ... V_m; U_0, ... U_n\}$, $Z$ must satisfy additional conditions. Loosely speaking, the key condition is that changes in $Y$ caused by varying $X$ should be uncorrelated with changes in $X$ caused by varying $Z$.

<u>Example 3</u> (continued): In this example, both $Z$ and $X$ are binary. Condition (i) says that the subject-level treatment effect $\Delta_{j0}Y(\omega) = U_1(\omega) - U_0(\omega)$ is uncorrelated with $\Delta_{10}I(X(\omega) = 1)$ which represents the change in actual treatment choice when the subject's assigned treatment were to change from 0 to 1. Under condition (i), the treatment choice for a subject $\omega$ can change under different treatment assignment, but this change should not be correlated with the subject-level treatment effect, as $\omega$ varies across the population. Whether such a condition is reasonable depends on the application. For example, suppose a weight loss center wants to compare two weight-loss programs by offering a substantial fee discount for one of them (randomly drawn) for each new customer. Then Z represents the treatment selected for the discount and $X$ represents the actual treatment chosen by the customer. In this case, it may be reasonable to assume that the change in treatment choice depends mostly on whether the subject is sensitive to the cost of the treatment and is largely uncorrelated with the difference in potential outcomes under the two treatments. There is a large literature on how to handle non-compliance in randomized trials -- a subject $\omega$ is non-compliant if $X(\omega) \neq Z(\omega)$. It is known that causal inference in this situation will require assumptions on the "compliance" variables $I(V_1(\omega)=1)$ and $I(V_0(\omega) = 1)$. For example, under the assumption $I(V_0(\omega) = 1)$ is equal to 0 almost surely, one can estimate a "treatment effect" defined as the average of the subject-level treatment effect $U_1(\omega) - U_0(\omega)$ over the subpopulation of subjects who choose treatment 1 (Heckman 1990). Also, if $I(V_1(\omega)=1)$ is almost always greater than (or almost always less than) $I(V_0(\omega) = 1)$, then one can estimate the average of the subject-level treatment effects over the subpopulation of subjects who can be induced to change treatment choice when the treatment assignment Z changes (Imbens and Angrist 1994). These assumptions are appropriate in the contexts of those works, but they identify causal effects based on alternative definitions different from the usual definition in the potential outcome framework, see Imbens 2014 for discussion of alternative definition of causal effects.

### 3) The case of discrete $X$ and continuous $Z$

In this section, we assume that $X \in \{x_0, \dots, x_n\}$, $Z \in [0,1]$.

Let $\{V(z), z \in [0,1]\}$ be a collection of variables such that for each $z$, $V(z)$ is a categorical variable whose distribution is given by $P(V = x_j | Z = z) = q_j(z), j = 1, \dots n$. Suppose $Y, X, Z$ are related by the following structural equations

$$Y = \sum_{j=0}^{n} I(X = x_j) U_j, \quad X = \int_0^1 \delta_Z(z) V(z) dz \qquad (2)$$

where $\delta_Z$ denotes the delta function centered on Z. Under (2), $X(\omega)$ is generated by first drawing $Z(\omega)$ from the distribution of Z, then $X(\omega) = V(Z(\omega))$. Since the relation between Y and X remains the same as that in section 2, we use the same definition for the causal effect of X on Y, i.e. $\theta = (\theta_1, \dots \theta_n)$ where $\theta_j = E(U_j - U_0)$. Note that no restriction has been placed on the joint distributions of the $U_j$'s and $V(z)$'s, thus equation (2) represents a general causal model relating $Y, Z$ and $Z$. To formulate our next result, let

$\Delta_{j0} Y = U_j - U_0 =$ change in Y caused by changing X from $x_0$ to $x_j$

$I(X) =$ the n-vector whose $j^{th}$ component is $I_j(X) = I(X = x_j)$

$\Delta_{z,\delta} I_j(X) = I(V(z + \delta) = x_j) - I(V(z) = x_j) =$ change in $I_j(X)$ when Z is increased by $\delta$

$\mu(z) = E(Y|Z = z)$

$b(z) = \frac{d}{dz} \mu(z)$

$a(z) = \frac{d}{dz} E(I(X)|Z = z) = \frac{d}{dz} q(z)$

<u>Theorem 2</u>: Suppose Z is independent of $\{V(z), z \in [0,1], U_0, \dots U_n\}$ and that $\mu, q_1, \dots q_n$ are differentiable in $z$. Assume also that

i) $\Delta_{j0} Y \perp \Delta_{z,\delta} I_j(X)$ for all $\delta$ small, $z \in [0,1], j = 1, \dots n$
ii) the set of functions $\{a_1(z), \dots a_n(z)\}$ are not linearly dependent on each other

then $\theta$ is identifiable via the equation $b(z) = \sum_1^n \theta_j a_j(z)$.

Based on this result, we can use Z as an instrument for the estimation of the causal effect of X and Y, as follows:

Step 1: Perform nonparametric regression of $X$ on $Z$ to estimate $q_j(z), j = 1,..n$. This should be done as a multiclass prediction problem with smoothness regularization on the prediction probabilities $q_j(z)$, so that the derivatives of are also given from the regression. Obtain $a(z)$ from $q(z)$ by differentiation..

Step 2: Perform nonparametric regression of $Y$ on $Z$ to estimate $\mu(z)$, again with regularization to ensure smoothness. Obtain $b(z)$ from $\mu(z)$ by differentiation.

Step 3: Fit the equation $b(z) = \sum_1^n \theta_j a_j(z)$ to obtain the estimate for $\theta$.

Example 2 (continued): In this example, the structural equations are

$$Y = \beta_0 + \beta_1 X + \eta, \quad X = I(X^* > 0), \quad X^* = \alpha_0 + \alpha_1 Z + \varepsilon$$

where $\eta$ and $\varepsilon$ are zero-mean errors independent of $Z$. Then,

$$q(z) = P(X = 1|Z = z) = P(\varepsilon > -\alpha_0 - \alpha_1 z)$$

$$\mu(z) = \beta_0 + \beta_1 q(z)$$

The conditions of Theorem 2 are satisfied if the distribution of $\varepsilon$ is absolutely continuous (i.e. it has density).

### 4) The case of continuous $X$ and $Z$

In this section, $X \in \mathbf{X}$, $Z \in \mathbf{Z}$, where $\mathbf{X}$ and $\mathbf{Z}$ are bounded rectangles in $R^n$ and $R^m$ respectively. We assume

$$Y = f(X, U), \quad X = g(Z, U), \qquad (3)$$

where the function $f$ is smooth in $x$ in the sense that $f(x, u)$ has first and second order derivatives in $x$, and that $f$ and its derivatives are uniformly bounded. Likewise, each component of $g(z, u) = (g_1(z, u), ..., g_n(z, u))$ is assumed to be smooth in $z$. As in the introduction, $U$ denotes all other variables that may affect $Y$ and $X$. This setting is quite general as the only condition on the structural equations $f$ and $g$ is smoothness with respect to $x$ and $z$ respectively, and no restriction has been placed on the distribution of $U$.

To formulate the concept of causal effects of $X$ on $Y$, let $h(x, u)$ be the gradient of $f(x, u)$ with respect to $x$, (i.e. $h(x, u) = \frac{\partial f}{\partial x}(x, u)$ ), and define $\frac{\partial Y}{\partial X}$ and $\frac{\partial Y}{\partial x}$ to be the random variables $h(X, U)$ and $h(x, U)$ respectively.

<u>Definition</u>: The causal effect of $X$ on $Y$ at $X = x$ is $\theta(x) = E(\frac{\partial Y}{\partial x}) = Eh(x, U)$

For each sample point $\omega \in \Omega$, $\frac{\partial Y}{\partial x}(\omega)$ tells us how $Y(\omega)$ would change if we hold all other variables fixed and varying only $X(\omega)$, when the current value of $X(\omega)$ is $x$. According to the definition, $\theta(x)$ is just the population average of this change in $Y(\omega)$. Like the subject-level treatment effect $d_i$ in example 3, here $\frac{\partial Y}{\partial x}(\omega)$ is a counterfactual variable useful in the formulation of the causal effect $\theta$ but is not directly observable.

We now discuss how to use $Z$ as an instrument for the inference of $\theta(x)$. Fist we consider the case of when $X$ and $Z$ are real valued. Similar to $\frac{\partial Y}{\partial X}$, we define $\frac{\partial X}{\partial Z} = k(Z, U)$ and $\frac{\partial X}{\partial z} = k(z, U)$ where $k(z, u) = \frac{\partial g}{\partial z}(z, u)$. Note that in this case all these random variables are real valued.

Also, define the functions

$$b(z) = \frac{d}{dz} E(Y|Z = z)$$

$$a(z) = \frac{d}{dz} E(X|Z = z)$$

Consider the following conditions:

  i)  $Z$ is independent of $U$
 ii)  $\frac{\partial Y}{\partial X} \perp \frac{\partial X}{\partial Z}$ conditional on $Z = z$
iii)  $a(z) \neq 0$

<u>Theorem 3 (scalar case)</u>: Suppose both $X$ and $Z$ are real valued and conditions (i), (ii), (iii) hold, then $\phi(z) = E\left(\frac{\partial Y}{\partial X}\Big|Z = z\right)$ is identifiable via the equation

$$a(z)\phi(z) = b(z)$$

According to theorem 3, for $Z$ to be used as an instrumental variable, it should affect $X$ in a nontrivial way (condition iii) and should be independent of all other causal variables that affect $Y$ and $X$ (condition i). Furthermore, conditional on $Z$, the change in $Y$ caused by varying $X$ must be uncorrelated to the change in $X$ caused by varying $Z$ (condition ii).

<u>Example 4</u> (continued): In this example $Y = s(X) + U_1$, $X = t(Z) + U_2$, all variables are real valued, and $Z$ is independent of $(U_1, U_2)$. It follows that $a(z) = t'(z)$, $\frac{\partial Y}{\partial X} = s'(X)$, $\frac{\partial X}{\partial Z} = t'(Z)$. Clearly, $\frac{\partial Y}{\partial X}$ and $\frac{\partial X}{\partial Z}$ are uncorrelated conditional on $Z = z$. Suppose $t'(z) \neq 0$ for almost all $z$, then by theorem 3, we have

$$E(s'(X)|Z = z) = E\left(\frac{\partial Y}{\partial X}\Big|Z = z\right) = \phi(z) = b(z)/a(z) \quad (4)$$

Under the further assumption that that the family of conditional distributions $D(X|Z = z)$ is rich enough so that the conditional expectations determine any function $r(X)$ uniquely, then the function $s'()$ can be obtained by solving the integral equation (4). Finally the causal effect of $X$ on $Y$ is identified as $\theta(x) = E(\frac{\partial Y}{\partial x}) = E(s'(x)) = s'(x)$. We note that a similar condition (i.e. the conditional distributions are rich enough) was first introduced into causal inference by Newey and Powell, 2003).

<u>Example 5</u>: Consider the model $Y = U_1 X + U_2$, $X = t(Z) + U_3$, $Z$ is independent of $(U_1, U_2, U_3)$, $U_3$ is independent of $U_1$. In this model, both $\frac{\partial Y}{\partial X} = U_1 = \frac{\partial Y}{\partial x}$, and the causal effect $\theta$ is constant in $x$. Conditions (i)-(iii) can be verified as in the last example, so we can use theorem

3 to obtain $\phi(z) = E(\frac{\partial Y}{\partial X}|Z = z)$. Finally, the causal effect of $X$ on $Y$ is identified as $\theta = E(\phi(Z))$.

Under the general model (3), we can use theorem 3 to find $\phi(z) = E(h(X,U)|Z = z)$. Unfortunately, the function $\phi(z)$ cannot determine the function $\theta(x) = E(h(x,U))$ without additional assumptions. In example 4, it is assumed that $\frac{\partial Y}{\partial X}$ depends only on $X$ but not on $U$ and that $\frac{\partial X}{\partial Z}$ depends only on $Z$ but not on $U$. In example 5, it is assumed that $\frac{\partial Y}{\partial X}$ depends on $U$ but not on $X$. At this time we do not know what is the weakest assumption under which $\theta(\,)$ can be obtained from $\phi(\,)$.

Finally, we extend theorem 3 to the case when $X$ and $Z$ are vector valued (i.e. $n \geq 1, m \geq 1$). In this case, $\partial Y/\partial X_j = h_j(X,U)$ where $h_j(x,u) = (\frac{\partial f}{\partial x_i})(x,u)$, and $\partial X_j/\partial Z_i = k_{ij}(Z,U)$ where $k_{ij}(z,u) = (\frac{\partial g_j}{\partial z_i})(z,u)$, Furthermore, the functions $b(z)$, $a(z)$ and $\phi(z)$ take value in $R^m$ and $R^{m \times n}$ and $R^n$ respectively.

Consider the following conditions

i*)   $Z$ is independent of $U$

ii*)  For all $i$ and $j$, $\partial Y/\partial X_j \perp \partial X_j/\partial Z_i$ conditional on $Z = z$

iii*) $rank(a(z)) \geq n$

Theorem 3 (vector case): Suppose conditions (i*), (ii*), (iii*) hold, then $\phi(z) = E\left(\frac{\partial Y}{\partial X}\middle|Z = z\right)$ is identifiable via the equation

$$a(z)\phi(z) = b(z)$$

5) **Proofs**

Proof of theorem 1:

From (1), we have

$$Y = \sum_{j=0}^{n}(\sum_{i=1}^{m} I(Z = z_i)I(V_i = x_j))U_j$$
$$= \sum_{i=1}^{m} I(Z = z_i) (\sum_{j=0}^{n} I((V_i = x_j)U_j).$$

Let $s = \{j, k\}$, then

$$b_s = E(\sum_{j=0}^{n} I(V_i = x_j)U_j) - E(\sum_{j=0}^{n} I(V_k = x_j)U_j)$$

Replacing $I(V_i = x_0)$ by $1 - \sum_{j=1}^{n} I(V_i = x_j)$, we get

$$b_s = E(\sum_{j=1}^n (I(V_i = x_j) - I(V_k = x_j))(U_j - U_0))$$
$$= \sum_{j=1}^n E(I(V_i = x_j) - I(V_k = x_j))E(U_j - U_0)$$
$$= \sum_{j=1}^n (P(X = x_j|Z = z_i) - P(X = x_j|Z = z_k))E(U_j - U_0)$$
$$= \sum_{j=1}^n a_{sj}\theta_j$$

Proof of theorem 2:

$$E(I_j(X)|Z = z) = P(X = x_j|Z = z) = P(V(z) = x_j|Z = z) = P(V(z) = x_j) = q_j(z).$$

$$E(Y|Z = z + \delta) - E(Y|Z = z)$$
$$= E(\sum_0^n I(X = x_j)U_j|Z = z + \delta) - E(\sum_0^n I(X = x_j)U_j|Z = z)$$
$$= E(\sum_0^n I(V(z + \delta) = x_j)U_j) - E(\sum_0^n I(V(z) = x_j)U_j)$$
$$= \sum_1^n E(I(V(z + \delta) = x_j) - I(V(z) = x_j))(U_j - U_0)$$
$$= \sum_1^n E(I(V(z + \delta) = x_j) - I(V(z) = x_j))E(U_j - U_0)$$
$$= \sum_1^n (q_j(z + \delta) - q_j(z))\theta_j$$

Proof of Theorem 3:

Writing $E(X|z)$ instead of $E(X|Z = z)$ and ignoring terms of order $\delta^2$, we have

$$E(X|z + \delta) - E(X|z) = E(g(z + \delta, U)|Z = z + \delta) - E(g(z, U)|Z = z)$$
$$= E(g(z + \delta, U) - g(z, U)) = E\left(\frac{\partial g}{\partial z}(z, U)\right)\delta$$

Thus $a(z) = E\left(\frac{\partial g}{\partial z}(z, U)\right) = E(\frac{\partial g}{\partial z}(Z, U)|Z = z) = E(\frac{\partial X}{\partial z}|Z=z)$

Similarly,

$$E(Y|z + \delta) - E(Y|z)$$
$$= E(f(g(z + \delta, U), U) - f(g(z, U), U))$$
$$= E(\frac{\partial f}{\partial x}(g(z, U), U)(g(z + \delta, U) - g(z, U)))$$
$$= E(\frac{\partial f}{\partial x}(g(z, U), U)\frac{\partial g}{\partial z}(z, U))\delta$$

Thus $b(z) = E\left(\frac{\partial f}{\partial x}(g(z, U), U)\frac{\partial g}{\partial z}(z, U)\right)$
$$= E\{\left(\frac{\partial f}{\partial x}(g(z, U), U)\frac{\partial g}{\partial z}(z, U)\right)|Z = z\}$$
$$= E\{\left(\frac{\partial f}{\partial x}(g(Z, U), U)\frac{\partial g}{\partial z}(Z, U)\right)|Z = z\}$$

$$= E\left(\frac{\partial Y}{\partial X}\frac{\partial X}{\partial z}\Big|Z=z\right) = E\left(\frac{\partial Y}{\partial X}\Big|Z=z\right)E\left(\frac{\partial X}{\partial z}\Big|Z=z\right) = \phi(z)a(z)$$

Proof of Theorem 3*: Similar to the proof of Theorem 3.

**Acknowledgment:** The authors thanks Peng Ding, Xiaotong Shen and Dylan Small for reading a preliminary draft and providing helpful comments. This work was supported by NSF grants DMS1811920 and DMS1952386.

**References**

Angrist, J. D., Imbens, G. W. and Rubin, D. B. (1996) Identification of causal effects using instrumental variables. *Journal of the American Statistical Association*, 91: 444-455.

Chernozhukov, V., Imbens, G. W., Newey, W. K. (2007) Instrumental variable estimation of nonseparable models. *Journal of Econometrics* 139: 4-14

Durbin, J. (1954) Errors in variables. *Review of the International Statistical Institute*, 22:23-32

Heckman, J. J. (1990) Varieties of selection biases. *American Economic Review, Papers and Proceedings*, 80:313-338

Heckman, J. J. and Robb, R. (1985) Alternative methods for evaluating the impact of interventions, in *Longitudinal Analysis of Labor Market Data*, edited by J Heckman and B. Singer. New York: Cambridge University Press.

Imbens, G.W. (2014) Instrumental variables: an econometrician's perspective. *Statistical Science*, 3: 323-358

Imbens, G. W. and Angrist, J. D. (1994) Identification and estimation of local average treatment effects. *Econometrica,* 62:467-475.

Imbens, G.W. and Newey, W.K. (2009) Identification and estimation of triangular simultaneous equations models without additivity. *Econometrica*, 77:1481-1512.

Kennedy, E.H., Lorch, S. and Small, D.S. (2019) Robust causal inference with continuous instruments using the local instrumental variable curve. *J. R. Statist. Soc.* B, 81, Part 1, 121-143.

Newey, W. K. and Powell, J. (2003) Instrumental variable estimation of nonparametric models. *Econometrica*, 71:1565-1578.

Rubin, D (1974) Estimating causal effects of treatments in randomized and non-randomized studies. *Journal of Educational Psychology*, 66: 688-701

Torgovitsky, A (2015) Identification of nonseparable models using instruments with small support. *Econometrica*, 83:1185-1197